 \documentclass[letterpaper, 10 pt, conference]{ieeeconf} 

\IEEEoverridecommandlockouts
\usepackage{amsmath, amssymb, amsfonts}
\usepackage{enumerate}
\usepackage{setspace}
\usepackage{graphicx}
\usepackage{color}
\usepackage{verbatim}
\usepackage{url}
\usepackage{algorithm}
\usepackage[noend]{algpseudocode}
\usepackage{float}
\usepackage{amsthm,thmtools}

\makeatletter
\def\BState{\State\hskip-\ALG@thistlm}
\makeatother

\def\urltilde{\kern -.15em\lower .7ex\hbox{\~{}}\kern .04em}
\def\urldot{\kern -.10em.\kern -.10em}
\def\urlhttp{http\kern -.10em\lower -.1ex\hbox{:}\kern -.12em\lower 0ex\hbox{/}\kern -.18em\lower 0ex\hbox{/}}

\declaretheorem[name={Definition}  ] {definition} 
\declaretheorem[name={Theorem}  ] {theorem}

\declaretheorem[name={Problem}  ] {problem} 

\declaretheorem[name={Claim}  ] {claim}

\declaretheorem[name={Corollary}  ] {corollary}

\declaretheorem[name={Proposition}  ] {proposition}

\declaretheorem[name={Example},qed={\lower-0.3ex\hbox{$\square$}} ] {example}

\newcommand{\notop}{\bar}

\newcommand{\IEEEQED}{\hfill{$\square$}}

\newcommand{\andop}{\wedge}

\newcommand{\orop}{\vee}

\newcommand{\N}{\mathcal N}

\newcommand{\I}{\mathcal I}

\begin{document}
\title{Output Selection and Observer Design for Boolean Control Networks: \\A Sub-Optimal Polynomial-Complexity
 Algorithm\thanks{
			This research is supported in part by a research 
		grant from the Israel Science Foundation~(ISF grant 410/15).}} 
		
\author{Eyal Weiss  and Michael Margaliot\thanks{The authors are with the School of Elec.
 Eng., Tel-Aviv University,  Israel~69978. Corresponding Author: michaelm@eng.tau.ac.il}}

\maketitle
 	 
\begin{abstract}
 Using a graph-theoretic approach,
we derive 
 a new sufficient condition for observability of  a Boolean control network~(BCN).
Based on this condition, we describe  two algorithms: 
the first 
selects a set of   nodes so that observing this set makes the~BCN observable. 
The  second algorithm builds an   observer for the observable~BCN. 
Both algorithms are sub-optimal, as they are based on a sufficient but not necessary condition for observability.
 Yet their  time-complexity   is  {linear} in the 
length of the description of the~BCN, rendering them
feasible  for large-scale networks. 
We   discuss 
how these results can be used to provide a sub-optimal yet polynomial-complexity algorithm for the 
 minimal observability problem in~BCNs. Some of
the theoretical results are demonstrated using a~BCN
  model of the core network regulating the mammalian cell cycle.
\end{abstract}
\section{Introduction}\label{sec:intro}

 Boolean networks~(BNs) are discrete-time dynamical systems with Boolean state-variables~(SVs) 
and Boolean update functions. 
BNs have found wide applications in modeling and analysis of
 dynamical systems. 
They have been used to capture the existence and directions of links in 
complex systems (see, e.g.,~\cite{shahrampour2015topology}),
to model social networks (see, e.g.~\cite{soc_nets_bool,soc_nets_plos}), 
and the spread of epidemics~\cite{bool_epidemcs}.
 
Recently, BNs are extensively used in systems biology (see, e.g.~\cite{doi:10.1093/bioinformatics/btw747,1478-3975-9-5-055001}). 
 A typical  example is modeling gene regulation networks using~BNs (see, e.g.~\cite{KAUFFMAN1969437,othmer}). 
Here
  the state of each gene, that may be   expressed or not, is modeled using a Boolean SV.
The interactions between the genes 
(e.g. through the effect of the proteins that they encode on the promoter regions of other genes) determine
the Boolean update function for each~SV.

BNs with (Boolean) control inputs are called Boolean control networks~(BCNs).
  Cheng and his colleagues~\cite{cheng2001semi} used
	the semi-tensor product~(STP) of matrices 
to represent the dynamics of a~BCN   in a form similar  to  that of a discrete-time
linear control system.
 This led to   the analysis of many
 control-theoretic  problems for~BCNs, including
controllability~\cite{cheng2009controllability,dima_cont}
observability~\cite{fornasini2013observability,cheng2009controllability},
disturbance decoupling~\cite{cheng2011disturbance},
stabilization by state feedback~\cite{li2013state}
and optimal control~\cite{fornasini2014optimal,zhang2017finite,dima_opt}.
 However, the STP representation is   
exponential in the number of Boolean SVs, and thus it cannot be used in general 
to develop efficient algorithms for addressing these problems.


In many   real-world 
systems  it is not possible to measure the state of all
 the~SVs.  
For instance, the function of a multipolar neuron 
may depend upon signals received from thousands of other interconnected neurons
(see, e.g.,~\cite{hawkins2016neurons}).


State observers  are needed to reconstruct the entire state of the system based on a time sequence of what can be measured, i.e. the 
system outputs.
An estimate of the entire state is useful in many applications, as it greatly assists in monitoring and control. For example, in
bio-reactors  
one may use sensors to measure variables such as dissolved oxygen, pH and temperature, yet
 key biological state-variables such as biomass and product concentrations can be much more difficult (and costly) to measure~\cite{hulhoven2006hybrid}. 
A typical application of state-observers is to integrate
 their estimates  in full-state feedback controllers (see, e.g.,~\cite{sontag_book}).

The existence of  observers is known to depend  on the reconstructability property of the system,
which was shown to be equivalent to another property, called observability,
 in linear systems~\cite{kalman1969topics}
and in BNs~\cite{fornasini2013observability}.
In BCNs, there are several different definitions for observability appearing in the
literature~\cite{cheng2009controllability,zhao2010input,dima_obs,fornasini2013observability}
(see also the work on   final
state observable graphs  in~\cite{obsgraphs}),
making the connection to reconstructability somewhat complex.
In subsection~\ref{ssec:Observability} we state the definition used in this paper and clarify 
its relation to the definition of reconstructability. 

When a given system is not observable, it is sometimes possible to make it observable by placing
additional sensors that measure
more (functions of the) SVs. Of course, this may be costly in terms of resources,
so a natural question is:
find the \emph{minimal} number of measurements to add so that the resulting system is observable. 
This \emph{minimal observability problem} is also interesting theoretically, as its
solution means identifying the (functions of) SVs that provide the maximal
information on the entire state of the system~\cite{Liu_observ}.
Indeed, minimal observability problems have recently attracted considerable interest. 
Examples include monitoring complex services by minimal logging~\cite{biswas2008minimal}, 
the optimal placement of phasor measurement units   
in power systems (see, e.g.,~\cite{peng2006optimal}), and the minimal sparse 
observability problem addressed in~\cite{sarma2014minimal}. 


It is known  that testing observability of BCNs is NP-hard 
in the number of~SVs of the system~\cite{dima_obs} 
(for a general survey on the computational complexity of various problems in 
systems and control theory, see~\cite{blondel}). 
This means that, unless P=NP,
it is computationally intractable to determine whether a large BCN is observable.
This   implies that the minimal observability problem in BCNs is also
NP-hard, since it must entail analyzing  observability.

It is thus not surprising that many observers for~BCNs have exponential complexity. 
These include  the Shift-Register observer, 
the Multiple States observer~\cite{fornasini2013observability}
and the Luenberger-like observer~\cite{zhang2016observer}.
This complexity implies that these algorithms cannot be used in 
  large-scale networks.

This paper is motivated by recent work on the minimal observability problem in a special class
 of~BNs
called \emph{conjunctive Boolean networks}~(CBNs). In a CBN every   update function 
 is comprised of only AND operations. 
A necessary and sufficient observability condition for~CBNs
has been derived in~\cite{weiss2017polynomial}. This condition is based on a graph representing  
the dependencies between the SVs and the update functions. 
This yields polynomial-time algorithms for:
 (1)~solving the minimal observability problem in~CBNs;
and (2)~designing observers for observable CBNs.

Here, we show that a similar graph-theoretic approach 
   provides a sufficient (but not necessary) 
 condition for observability of (general)~BCNs. 
This induces algorithms for solving Problems~(1) and~(2) above.
Now the algorithms are not optimal anymore, i.e. they may add  more observations than the minimal number that is indeed
required. Nevertheless,
they retain their polynomial-time complexity implying that they are feasible for large-scale BCNs.
 Of course, in the particular case where all the update functions 
are~AND gates  these algorithms become optimal.

We note in passing that the special structure of CBNs makes them amenable to analysis
 (see, e.g. ~\cite{basar_CBN,weakly_con_cbn,basar_min_cont,min_control})
and we believe that more results from this field can be extended to handle general~BNs. 
 
The next section reviews some known definitions and results that are used
in deriving the main results in Section~\ref{sec:main}.
 
\section{Preliminaries}\label{sec:pre}

\subsection{Boolean Control Networks}\label{ssec:BCN}
 
Let~$S:=\{0,1\}$. For two integers~$i,j$ let~$[i,j]:=\{i,i+1,\dots,j\}$. 
A BCN with $n$ SVs, $p$ inputs and~$m$ outputs can be represented by the following equations:
\begin{align}\label{eq:BCN}
X_i(k+1)&=f_i(X(k),U(k)),  \quad &\forall i\in[1,n],\nonumber\\
Y_j(k)&=h_j(X(k)), \quad &\forall j\in[1,m], 
\end{align}
where the state vector at time~$k$ is denoted
by~$X(k):=\begin{bmatrix} X_1(k)&\dots&X_n(k)\end{bmatrix}' \in S^n$,
the input vector by ~$U(k):=\begin{bmatrix} U_1(k)&\dots&U_p(k)\end{bmatrix}' \in S^p$,
the output vector by~$Y(k):=\begin{bmatrix} Y_1(k)&\dots &Y_m(k)\end{bmatrix}' \in S^m$
and $f_i$, $h_j$ are Boolean functions.

In principle, an update function may include an SV (or control input) 
that has no effect on the function, e.g.~$f_1(X_1(k),X_2(k))=X_1(k) \orop ( X_2(k) \andop  \notop  {X}_2(k)  )$.
We assume that    such arguments have been removed, in other words, if
$X_i(k)$ 
is an argument of $f_j$ then  there exists an assignment 
 of~$X_1,\dots,X_{i-1},X_{i+1},\dots ,X_n,U$
such that
\begin{align*}
	f_j(X_1,\dots,X_{i-1},0,X_{i+1} ,\dots ,X_n ,U ) \\
	\neq f_j(X_1,\dots,X_{i-1} ,1,X_{i+1} ,\dots,X_n ,U ).  
\end{align*}
Similarly, a control input appears in an update function only if it really  affects the function.

If for some~$i$ there exists an output  
$Y_j(k)=X_i(k)$ then we say that the SV~$X_i$ 
is \emph{directly observable} or \emph{directly measurable}.

\subsection{Observability and   Reconstructability}\label{ssec:Observability}
Observability [reconstructability] refers to the ability to uniquely determine
the   {initial} [final] state of a system based
on a time-sequence of the input and output.
There are several definitions in the literature for 
observability of BCNs.
 Here, we follow the one used in~\cite{fornasini2013observability}.

 \begin{definition}\label{def:observability} 
 	We say that~\eqref{eq:BCN} is \emph{observable} on~$[0,N]$ 
 	if for every control sequence $  U:=\{  U(0),\dots,  U(N-1)\}$ and
 	for any two different initial conditions~$X(0)$ and~$\tilde X(0)$
 	the corresponding solutions of the BCN for~$  U$ yield \emph{different} 
 	output sequences~$\{Y(0),\dots,Y(N)\}$ and~$\{\tilde Y(0),\dots, \tilde Y(N) \}$.
 \end{definition}

 This means that 
  it is always possible 
to uniquely determine the initial condition  from the output sequence on~$[0,N]$.
A~BCN is called \emph{observable} if it is observable for some value~$N \geq 0$.


 \begin{definition}\label{def:reconstructability}~\cite{fornasini2013observability}
	The BCN~\eqref{eq:BCN} is said to be reconstructable 
	if for any integer~$r>0$   the knowledge of every admissible input and output
	trajectories~$\{(Y(k),U(k))\}$, $ k = 0,1,\dots,r$, uniquely determine
	the final state $X(r)$.
 \end{definition}
 Clearly, if a BCN is observable then it is also reconstructable.

\subsection{Directed Graphs}\label{ssec:digraph}

Let $G=(V,E)$ be a directed graph (digraph), with~$V$ the set of vertices, 
and $E$ the set of directed edges (arcs). Let~$e_{i\to j}$ (or $(v_i\to v_j)$) denote
the arc from~$v_i$ to~$v_j$. When such an arc exists, we say that~$v_i$ is an \emph{in-neighbor} of~$v_j$, and~$v_j$ as 
an \emph{out-neighbor} of~$v_i$. 
The set of in-neighbors [out-neighbors] of~$v_i$ is denoted by~$\N_{in}(v_i)$ [$\N_{out}(v_i)$]. 
The \emph{in-degree} [\emph{out-degree}] 
of~$v_i$ is~$|\N_{in}(v_i)|$ [$|\N_{out}(v_i)|$].
A \emph{source} [\emph{sink}] is a node with  in-degree [out-degree] zero.

For  $v_i,v_j\in V$ a \emph{walk} from~$v_i$ to~$v_j$, denoted~$w_{ij}$, is a sequence: 
$v_{i_0}v_{i_1}\dots v_{i_q}$, with $v_{i_0}=v_i$, $v_{i_q}=v_j$,  
and~$e_{i_k\to i_{k+1}} \in E $ for all $k\in [0,q-1]$. 
A \emph{closed walk} is a walk that starts and terminates at the same vertex. 
A closed walk is called a \emph{cycle} if all the vertices in the walk are distinct, except for the start-vertex and the end-vertex.
 
\subsection{Dependency Graph}\label{ssec:dependency_graph}

The dependency graph of a BCN is defined by $G=(V,E)$,
where $V=\{X_1,\dots,X_n,U_1,\dots,U_p\}$ i.e. every vertex corresponds to either an SV or an input of the BCN.
 The edge
$(X_i \to X_j)\in E$ [$(U_q \to X_j)\in E$] iff $X_i(k)$ [$U_q(k)$]
is an argument of $f_j$ (the update function of $X_j(k+1)$).
 Thus, the dependency graph encodes the actual variable dependencies in the update functions.

We denote the subset of nodes corresponding to SVs by~$V_s:=\{X_1,\dots,X_n\} \subseteq V$,
the subset of edges involving only~SVs by~$E_s \subseteq E$, and let~$G_s:=(V_s,E_s)$
be the  resulting ``reduced'' dependency graph.
Referring to the outputs of the~BCN, we say that a node in the dependency graph that
represents a [non] directly observable SV    
is   a \emph{[non] directly observable node}.

\section{Main Results}\label{sec:main}

From here on, we consider BCNs   in the form: 
\begin{align}\label{eq:BCN_regular_form}
X_i(k+1)&=f_i(X(k),U(k)),  \quad  &i\in \{1,\dots,n\},\nonumber\\
Y_j(k)&=X_j(k), \quad &j\in\{1,\dots,m\}, 
\end{align}
  that is, every output~$Y_j$ is the value of an~SV. 
We assume without loss of generality 
that the~$m$ outputs correspond to the first~$m$ SVs. 
Thus, nodes~$X_1,\dots,X_m$ [$X_{m+1},\dots,X_n$] in the dependency graph
are [non] directly observable. 

\subsection{Sufficient Condition for Observability}\label{ssec:conditions}

We begin by presenting two definitions which will be used later on.

\begin{definition}\label{def:P_1}
	We say that a BCN has \emph{Property~$P_1$} if for  every non-directly observable node~$X_i$ 
	there exists some other node~$X_j$  such that
	$\N_{in}(X_j)=\{X_i\}$.
\end{definition}
In this case,~$X_j(k+1)=f_j(X_i(k))$, and thus
  either~$X_j(k+1)= X_i(k)$ or~$X_j(k+1)=\notop X_i(k)$. 
 This means that  the information on the state of~$X_i$ 
  ``propagates'' to~$X_j$.

\begin{definition}\label{def:P_2}
	We say that a BCN has \emph{Property~$P_2$} if
		every cycle~$C$ in its dependency graph     
	that is composed solely of non-directly observable nodes 
  satisfies the following property:
	$C$ includes a node~$X_i$ which is the only element in the in-neighbors set of some other
	node~$X_j$, i.e.~$N_{in}(X_j)=\{X_i\}$, and~$X_j$ is not part of the cycle~$C$.
\end{definition}
This means that~$X_j(k+1)=f_j(X_i(k))$, so  the information on the state of every node in the cycle
propagates to~$X_i$
 and then  to~$X_j$, where~$X_j$ is not part of the cycle.

We now provide a     sufficient condition for observability.
\begin{theorem}\label{theorem:sufficient_condtion}
	A BCN that satisfies properties~$P_1$ and~$P_2$ is observable.
\end{theorem}

\noindent To prove this, we first introduce another definition and several auxiliary results. 

\begin{definition}\label{def:obsereved_paths}
	An \emph{observed  path} in the dependency graph is a non-empty ordered set of nodes  
	such that: (1)~the last element in the set is a directly observable node; and (2)~if 
	the set contains~$p>1$ elements, then for any~$i<p$ 
	the~$i$-th element is a non-directly observable node, 
	and is the \emph{only} element in the in-neighbors set of node~$i+1$.
	Observed paths with non-overlapping nodes are called \emph{disjoint observed paths}.	
\end{definition}

\noindent Roughly speaking, an observed path corresponds to a ``shift register''
whose last cell is directly observable.

\begin{proposition}\label{prop:vertex_cover_by_OPs}
	Consider a BCN   that satisfies properties~$P_1$ and~$P_2$. Then:
	\begin{enumerate}
		\item $G_s$ can be decomposed  
		into disjoint observed paths, such that every vertex~in $G_s$
		belongs to a single observed path 
		(i.e., the union of the disjoint observed paths is a vertex cover of~$G_s$).
		\item For every vertex~$v \in (V \setminus V_s)$, $\N_{out}(v)$ contains only
		vertices that  are located at the beginning of observed paths.	  
	\end{enumerate}	
\end{proposition}

{\sl Proof of Prop.~\ref{prop:vertex_cover_by_OPs}}.
We give a constructive proof.  Algorithm~1   below 
accepts   a graph~$G$ that satisfies properties~$P_1$ and~$P_2$
and terminates after each vertex in~$G_s$ 
belongs to exactly one observed path. Comments in the algorithm are enclosed within~(* $\ldots$ *).

\begin{algorithm}[H] \label{alg:one} 
\caption{Decompose the nodes of~$G_s$ into disjoint observed paths} 
\renewcommand{\algorithmicrequire}{\textbf{Input:}}
\renewcommand{\algorithmicensure}{\textbf{Output:}}
\begin{algorithmic}[1]
\Require Dependency graph~$G$ of a BCN in the form~\eqref{eq:BCN_regular_form}
 that satisfies properties~$P_1$ and~$P_2$.
\Ensure A decomposition of~$G_s$ into $m$ disjoint observed paths.
\For  {$i=1$ to $m$}
  (* every iteration  builds a new path ending with~$X_i$ *) 
\State  $\textit{o-node} \gets X_i$ ;
   $\textit{o-path}\gets \{X_i\}$
\If {$|\N_{in}(\textit{o-node})|=1$} \label{if:in_deg_equals_one}
\State Let $v$ be   such that~$\{v\}=\N_{in}(\textit{o-node})$
	\If {$v  $ does not belong to a  previous path, $v \in V_s$ \\\hspace*{1.5cm} and~$v$ is not 
		 directly observable}  
		\State   insert~$v$ to $\textit{o-path}$ just before~$\textit{o-node}$
		\State $\textit{o-node} \gets v$;
				  goto~\ref{if:in_deg_equals_one}
	\Else { print} $\textit{o-path}$
	\EndIf 
\EndIf
\EndFor
\State  {\bf end for}
\end{algorithmic}
\end{algorithm}

We now prove the correctness of Algorithm~1. 
To simplify the notation, we say that~$v_p$ \emph{points to}~$v_q$ if~$p\not =q$ and~$N_{in}(v_q)=\{v_p\}$, and
  denote this by~$v_p \mapsto v_q$. 
  The special arrow indicates that the dependency graph includes an edge from~$v_p$ to~$v_q$
	and that there are no other edges pointing to~$v_q$. 

If all the nodes of~$G_s$ are directly observable (i.e. if~$m=n$) 
  the algorithm will assign every node to a different observed path and this is correct.
Thus, we may assume that~$m<n$. 
Pick   a non directly observable node~$X_j$. Then~$m<j\leq n$. 
Our first goal is to  prove the  following result.
\begin{claim}\label{cla:out}
 The algorithm outputs  an observed path  that contains~$X_j$. 
\end{claim}
By Property~$P_1$, there exists~$k\not = j$ such that  $X_j \mapsto X_k$. 
  We consider two cases. 

\noindent {\sl Case 1.} 
If~$k\leq m$ then~$X_k$ is directly observable and
  the algorithm will add~$X_j$ to an observed path as  
	it ``traces back'' from~$X_k$ unless~$X_j$ has already been included in some other observed path found by the algorithm. 
	Thus, in this case Claim~\ref{cla:out} holds. 

\noindent {\sl Case 2.}  Suppose that~$k>m$, i.e.~$X_k$ is non directly observable.
By Property~$P_1$, there exists~$h\not = k$ such that~$ X_k \mapsto X_h $, so~$X_j \mapsto X_k  \mapsto X_h $.  
If~$h\leq m$ then we conclude as in Case~1 that 
 the algorithm outputs  an observed path  that contains~$X_j$. 
Thus, we only need to consider the case where as we proceed from~$X_j$ using Property~$P_1$ we never ``find'' a
directly observable node. Then 
there exists a set of non directly observable nodes~$X_{k_1},\dots, X_{k_\ell}$, with~$k_1=j$,
 such that 
$
													  X_{k_1} \mapsto X_{k_2} \mapsto \dots \mapsto 
														X_{k_\ell} 	\mapsto	   X_{k_1}  .
$
This means that~$X_j$ is   part of a cycle~$C$ of non directly observable nodes. 
 By Property~$P_2$, $C$ includes  a node~$X_{k_{i}}$ such that~$ X_{k_{i}} \mapsto X_{s_1}  $,
  where~$X_{s_1}$ is not part of the cycle~$C$.
If~$X_{s_1}$ is   directly observable then we conclude that the algorithm will output 
an   observed path that includes~$X_j$. 
If~$X_{s_1}$ is not directly observable then by  Property~$P_1$,  there exists~$s_2 \not =s_1 $ such that~$X_{s_1} \mapsto X_{s_2}$. Furthermore, since every node in~$C$ has in degree one,~$X_{s_2} \not \in C $. Proceeding this way, 
we conclude that there exist~$s_1,\dots , s_p$ such that
$
					X_{k_{i}} \mapsto X_{s_1}  \mapsto X_{s_2}\mapsto  \dots\mapsto X_{s_p},
$
with~$X_{s_p}$ a directly observable node. This means that the algorithm will output~$X_j$
in an observed path  as it traces back from~$X_{s_p}$, unless it already included~$X_j$  
 in another observed path. 
This completes the proof of  Claim~\ref{cla:out}.~\IEEEQED

Summarizing, we showed that  \emph{every} non directly observable node~$X_j$ is contained in an observed path produced by the algorithm. 
The fact that every~$X_j$ will be in a single observed path, and that the observed paths will be distinct is clear from the description of the algorithm.
From the definition of an observed path, it is clear that only a vertex 
which is located at the beginning of an observed path may contain edges from 
vertices in $V \setminus V_s$. 
This completes the proof of Prop.~\ref{prop:vertex_cover_by_OPs}.~\IEEEQED  

We can now prove Thm.~\ref{theorem:sufficient_condtion}. 
 
{\sl Proof of Thm.~\ref{theorem:sufficient_condtion}.}
Consider   the  following three    statements:
\begin{enumerate}[(a)]
\item The dependency graph~$G$ has Properties~$P_1$ and~$P_2$;
\item There   exists a decomposition 
of the dependency graph  
into~$m\geq 1$ disjoint observed paths~$O^1,\dots,O^m$, 
such that every vertex in~$G$ belongs to a single observed path.
\item The BCN is observable.
\end{enumerate}

The  correctness of Algorithm~1 
implies that~$(a) \to (b)$.  We now show that~$(b) \to(c)$.
Suppose that~$(b)$ holds. Let~$O^i=(X_{i_1},\dots,X_{i_{N_i}})$ be an observed path. 
The output  of  this path at times~$0, \dots, N_i-1$,  is
\begin{align*}
X_{i_{N_i}}(0)&=X_{i_{N_i}}(0),\\
 X_{i_{N_i}}(1)&=f_{i_{N_i}}(X_{i_{N_i-1}}(0)),\\
X_{i_{N_i}}(2)&=f_{i_{N_i}}(f_{i_{N_i-1}}(X_{i_{N_i-2}}(0))),\\ & \vdots\\
X_{i_{N_i}}(N_i-1)&=f_{i_{N_i}}(f_{i_{N_i-1}}(\dots f_{i_2}(X_{i_1}(0)))).
\end{align*}
Each of $f_{i_1},\dots,f_{i_{N_i}}$ is either the identity function or 
the NOT operator,  so  they are all invertible.
Moreover, their inverse functions are the functions themselves,
namely, $f_{i_1}^{-1}=f_{i_1},\dots,f_{i_{N_i}}^{-1}=f_{i_{N_i}}$.
Hence, the initial values of the~SVs in~$O^i$ can be reconstructed
as~$X_{i_{N_i}}(0)$,
$X_{i_{N_i-1}}(0)=f_{i_{N_i}}(X_{i_{N_i}}(1))$,
$X_{i_{N_i-2}}(0)=f_{i_{N_i-1}}(f_{i_{N_i}}(X_{i_{N_i}}(2)))$, and so on.
We conclude that 
  it is possible to determine the initial condition of every~SV
in  the~BCN using the output sequence on~$[0,\max_{i=1,\dots,m} \{N_i \} -1]$,
and this holds for \emph{every} control sequence.
Thus, the~BCN is observable, so $(b) \to (c)$. We conclude that~$(a) \to (c)$
 and this proves
Thm.~\ref{theorem:sufficient_condtion}.~\IEEEQED

\subsection{Observer Design}\label{ssec:dpatho}

The proof of Thm.~\ref{theorem:sufficient_condtion} implies the following result. 
\begin{corollary}\label{coro:observer}
Consider  a BCN that satisfies the sufficient condition stated in Thm.~\ref{theorem:sufficient_condtion}.
An  observer for this BCN can be designed as follows:
\begin{enumerate} [(a)]
	\item Construct the dependency graph~$G$;
	\item Apply Algorithm~$1$ to decompose the nodes of~$G_s$ into a set of disjoint observed paths;
	\item Observe an output sequence of length equal to the longest observed path;
	\item Map the values observed at each output to the values of the SVs composing the observed paths, 
	as done in the proof of Thm.~\ref{theorem:sufficient_condtion}, 
	to obtain the initial state~$X(0)$.
\end{enumerate}
\end{corollary}

Of course, once the initial condition is determined the state for all time~$k$ can be obtained using the known input sequence and dynamics. 
We refer to the observer described above
   as the \emph{Disjoint-Path Observer}.

\begin{figure}[t]
	\centering
	
\begingroup%
\makeatletter%
\providecommand\color[2][]{%
	\errmessage{(Inkscape) Color is used for the text in Inkscape, but the package 'color.sty' is not loaded}%
	\renewcommand\color[2][]{}%
}%
\providecommand\transparent[1]{%
	\errmessage{(Inkscape) Transparency is used (non-zero) for the text in Inkscape, but the package 'transparent.sty' is not loaded}%
	\renewcommand\transparent[1]{}%
}%
\providecommand\rotatebox[2]{#2}%
\newcommand*\fsize{\dimexpr\f@size pt\relax}%
\newcommand*\lineheight[1]{\fontsize{\fsize}{#1\fsize}\selectfont}%
\ifx\svgwidth\undefined%
\setlength{\unitlength}{162.71521861bp}%
\ifx\svgscale\undefined%
\relax%
\else%
\setlength{\unitlength}{\unitlength * \real{\svgscale}}%
\fi%
\else%
\setlength{\unitlength}{\svgwidth}%
\fi%
\global\let\svgwidth\undefined%
\global\let\svgscale\undefined%
\makeatother%
\begin{picture}(1,0.92855473)%
\lineheight{1}%
\setlength\tabcolsep{0pt}%
\put(0,0){\includegraphics[width=\unitlength]{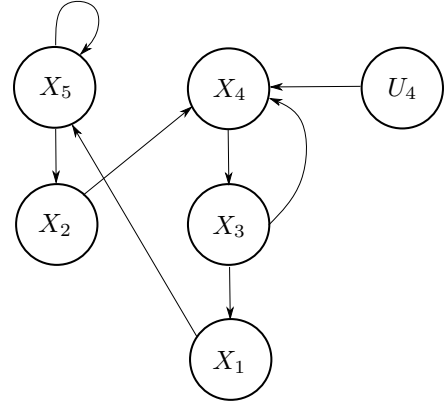}}%
\put(-0.21368419,0.98277574){\color[rgb]{0,0,0}\makebox(0,0)[lt]{\begin{minipage}{0.29631081\unitlength}\raggedright \end{minipage}}}%
\put(-0.82605993,2.15484974){\color[rgb]{0,0,0}\makebox(0,0)[lt]{\begin{minipage}{0.79016232\unitlength}\raggedright \end{minipage}}}%
\put(-0.0556517,1.93097044){\color[rgb]{0,0,0}\makebox(0,0)[lt]{\begin{minipage}{0.0395081\unitlength}\raggedright \end{minipage}}}%
\put(0.48429243,0.12018195){\color[rgb]{0,0,0}\makebox(0,0)[lt]{\begin{minipage}{1.14573536\unitlength}\raggedright \end{minipage}}}%
\put(0.05858424,0.39712341){\color[rgb]{0,0,0}\makebox(0,0)[lt]{\lineheight{1.25}\smash{\begin{tabular}[t]{l}$X_2$\end{tabular}}}}%
\put(0.46291507,0.71140635){\color[rgb]{0,0,0}\makebox(0,0)[lt]{\lineheight{1.25}\smash{\begin{tabular}[t]{l}$X_4$\end{tabular}}}}%
\put(0.06079407,0.71557588){\color[rgb]{0,0,0}\makebox(0,0)[lt]{\lineheight{1.25}\smash{\begin{tabular}[t]{l}$X_5$\end{tabular}}}}%
\put(0.46245833,0.08032941){\color[rgb]{0,0,0}\makebox(0,0)[lt]{\lineheight{1.25}\smash{\begin{tabular}[t]{l}$X_1$\end{tabular}}}}%
\put(0.4618948,0.3936191){\color[rgb]{0,0,0}\makebox(0,0)[lt]{\lineheight{1.25}\smash{\begin{tabular}[t]{l}$X_3$\end{tabular}}}}%
\put(0.86692566,0.71381892){\color[rgb]{0,0,0}\makebox(0,0)[lt]{\lineheight{1.25}\smash{\begin{tabular}[t]{l}$U_4$\end{tabular}}}}%
\end{picture}%
\endgroup%
	
	\caption{The dependency graph $G$ of the BCN in Example~\ref{ex:alg_one}.}
	\label{fig:alg_one}
\end{figure}

\begin{example}\label{ex:alg_one}
	Consider the single-input and two-output BCN:
	\begin{align}
	X_1(k+1)&=X_3(k),\nonumber \\
	X_2(k+1)&=\notop X_5(k),\nonumber\\
	X_3(k+1)&=\notop X_4(k),\nonumber\\
	X_4(k+1)&=(X_2(k) \andop X_3(k)) \orop U_4(k),\\
	X_5(k+1)&= \notop X_1(k)  \orop X_5(k),\nonumber\\
	Y_1(k) &=X_1(k),\nonumber\\
	Y_2(k) &=X_2(k)\nonumber.
	\end{align}
	The dependency graph~$G$  of this BCN  is  depicted  Fig.~\ref{fig:alg_one}, and it is straightforward to verify
	that it 
	satisfies Properties~$P_1,P_2$.
	 Thm.~\ref{theorem:sufficient_condtion} implies that this BCN   is observable,
	and that $G_s$ is decomposable to a set of disjoint observed paths.
	Applying Algorithm~$1$ to this BCN yields a decomposition to two observed paths:~$O^1=(X_4,X_3,X_1)$,
	$O^2=(X_5,X_2)$, where~$X_4 \mapsto X_3 \mapsto X_1$,
	$X_5 \mapsto X_2$. 
\end{example}

\subsubsection*{Complexity Analysis of the Construction of a Disjoint-Path Observer}
The complexity of generating the dependency graph~$G$ is  
linear in the description of the graph, which is~$O(|V|+|E|)$.
The resulting graph satisfies~$|V|=n+p$, $|E| \leq n^2+pn$,
with~$n$ being the number of SVs and $p$ the number of control inputs.
By using a chained-list data structure for representing the dependency graph,
Algorithm~$1$ can be implemented in~$O(n)$ time.
Indeed Algorithm~$1$ systematically passes through each vertex of~$V_s$ once
(with $|V_s|=n$) and performs~$O(1)$ operations on each such vertex.
After the decomposition to disjoint observed paths takes place,
determining the initial condition of the~BCN is attained in
complexity which is at most linear in the
length of the description of the~BCN.
From this stage, determination of the state at following time steps is
obtained by direct calculation of the dynamics, which is again
linear in length of the description of the~BCN for every time step.

Summarizing, the complexity of the construction is linear in the
length of the description of the~BCN, namely~$O(|V|+|E|)$. Thus, it satisfies the 
bound~$O(n^2+pn)$.

\subsection{Output Selection and an Upper Bound for Minimal Observability}\label{ssec:min_obs}
In some cases, it is possible to add sensors to measure more~SVs.
We consider the following problem.

\begin{problem}\label{prob:minobs}
	Given a BCN with~$n$ SVs determine a \emph{minimal} set of indices~$\I \subseteq [1,n]$,
	such that making each~$X_i(k)$, $i \in \I$, 
	directly measurable yields an observable~BCN.
\end{problem}

As mentioned in Section~\ref{sec:intro}, this minimal observability problem in BCNs is NP-hard.
We use the conditions in Thm.~\ref{theorem:sufficient_condtion}  
to provide  a sub-optimal yet nontrivial solution to Problem~\ref{prob:minobs}.
This is described by  
  Algorithm~$2$ which is discussed below.
 Furthermore, due to Corollary~\ref{coro:observer},
an efficient (polynomial-time) observer can also be obtained for a BCN 
with outputs chosen according to this scheme.

Algorithm~$2$ below is based on
an algorithm given in~\cite{weiss2017polynomial} that
 solves
the minimal observability problem in a special case of BCNs.
In our case (general BCNs) it does not offer a minimal solution to Problem~\ref{prob:minobs},
but rather it solely gives a solution 
which satisfies the conditions in Thm.~\ref{theorem:sufficient_condtion}.

The idea behind the algorithm is simple: it first creates three lists. A list~$L_1$   of all~SVs that are 
not directly observable and are not  the only element in the in-neighbors' set of another   node;
  a list~$L_2$  of all   SVs that are 
not directly observable and   are the  only element in the in-neighbors  set of another   node;
and a list~$L_C$ of   cycles composed solely out of nodes in~$L_2$. For each cycle~$C\in L_C$, it then
checks if one of its elements appears as the
 only element in the in-neighbors  set of another   node that is not part of~$C$.
 If so, it removes~$C$ from~$L_C$. Finally, it returns the~SVs in~$L_1$ 
and  one  element from each cycle~$C\in L_C$. Making these~SVs directly observable, 
by adding them as outputs, yields an observable~BCN.

Since the steps of the algorithm are basically described in~\cite{weiss2017polynomial},
yet here it is used for general~BCNs, we describe only the functional change 
that  is depicted in the headline (i.e., input-output description).

\begin{algorithm}[H] \label{alg:two} 
\caption{Output selection for meeting the condition of 
	Thm.~\ref{theorem:sufficient_condtion}: a high-level description}  
\renewcommand{\algorithmicrequire}{\textbf{Input:}}
\renewcommand{\algorithmicensure}{\textbf{Output:}}
\begin{algorithmic}[1]
\Require A BCN~\eqref{eq:BCN_regular_form} with~$n$ SVs and~$m\geq 0$ outputs.
\Ensure A set of SVs so that
		making these SVs directly observable yields a BCN that  
		satisfies conditions~$P_1$ and~$P_2$. 
\end{algorithmic}
\end{algorithm}

The complexity analysis of Algorithm~$2$ is done in~\cite{weiss2017polynomial},
yielding a runtime which is linear in the
length of the description of the~BCN, namely~$O(|V|+|E|)$.
We note that the algorithm provides a specific solution which meets the conditions of 
Thm.~\ref{theorem:sufficient_condtion}, 
but it is straightforward to modify this so that the algorithm will return 
the information needed to build various possible solutions which meet the conditions of theorem
(this is explained in detail in~\cite{weiss2017polynomial}). 
If the algorithm  returns an output list that is empty then the~BCN is observable, 
so it can also be used to determine if a given~BCN  
  satisfies the sufficient condition for observability.

%
%

\begin{theorem}\label{theorem:alg_validity}
  Algorithm~$2$ provides a solution that  satisfies the conditions of 
  Thm.~\ref{theorem:sufficient_condtion}.
\end{theorem}
The proof of Thm.~\ref{theorem:alg_validity} is the same one that appears in~\cite{weiss2017polynomial}
(but this time it has a different meaning, namely, here  we consider   the general class of BCNs),
so we omit it.
The correctness of Algorithm~$2$ implies the following:
\begin{corollary}\label{coro:upper_bound}
An upper bound for the minimal size of the solution to Problem~\ref{prob:minobs}
is given by the size of the solution generated by Algorithm~$2$.
\end{corollary}
Note that for the particular case of CBNs  the  
upper bound  
provided by size of the solution generated by Algorithm~$2$ is tight. 
 \subsection{A Biological Example}
Faure et al.~\cite{faure06,faure_2009} derived 
a BCN model for the mammalian cell cycle  
 that includes nine SVs
representing the activity/inactivity  
 of nine different proteins and a single
  input corresponding to the activation/inactivation of a regulating protein 
	in the cell.
		The BCN dynamics is given by:
 \begin{align}\label{eq:boolean_model}
                    x_1(t+1)&= (  \bar{u} (t) \wedge  \bar x_3 (t) \wedge  \bar x_4 (t)     \wedge \bar x_9 (t)) \\
										& \;\; \vee(x_5(t) \wedge  \bar{u} (t) \wedge  \bar x_9 (t)), \nonumber \\
                    x_2(t+1)&=(\bar x_1 (t)  \wedge \bar x_4 (t) \wedge \bar x_9 (t))\vee(x_5(t) \wedge \bar x_1(t) \wedge \bar x_9 (t)),\nonumber \\
                    x_3(t+1)&=x_2(t)  \wedge \bar x_1(t), \nonumber\\
                    x_4(t+1)&=(x_2(t) \wedge \bar x_1(t) \wedge \bar x_6 (t) \wedge ( \overline{x_7(t)  \wedge x_8(t)}))  \nonumber \\ &\;\; \vee (x_4(t)  \wedge \bar x_1(t) \wedge \bar x_6
                    (t) \wedge ( \overline{x_7(t)  \wedge x_8(t)} )), \nonumber\\
                    x_5(t+1)&=(\bar{u}(t) \wedge \bar x_3 (t) \wedge \bar x_4 (t) \wedge \bar x_9 (t)) \nonumber \\
                       & \;\;\vee(x_5(t) \wedge ( \overline{x_3(t)  \wedge x_4(t)} ) \wedge   \bar{u}(t) \wedge \bar x_9 (t)), \nonumber \\
                    x_6(t+1)&=x_9(t), \nonumber\\
                    x_7(t+1)&=(\bar x_4 (t) \wedge \bar x_9 (t)) \vee  x_6(t) \vee (x_5(t) \wedge \bar x_9 (t)), \nonumber\\
                    x_8(t+1)&=\bar x_7 (t)\vee(x_7(t) \wedge x_8(t) \wedge (x_6(t) \vee x_4(t) \vee x_9(t))), \nonumber\\
                    x_9(t+1)&=\bar x_6 (t) \wedge \bar x_7 (t) . \nonumber
\end{align}

The input    is considered constant  in~\cite{faure06}, that is, 
either~$u(t)\equiv \text{True}$ or~$u(t)\equiv \text{False}$. Under this assumption, the~BCN
becomes two~BNs. The simulations in~\cite{faure06} show that when~$ u(t)\equiv \text{True}$  
 the corresponding~BN  admits a globally attracting periodic
trajectory  composed of~$7$ states. The sequence of state transitions along this trajectory qualitatively matches cell cycle progression.
For  $u(t)\equiv \text{False}$ the BN admits a single state that is globally attracting. This state corresponds to the G0 phase (cell quiescence).

Ref.~\cite{dima_obs} used the STP representation 
 combined with a trial and error approach to conclude that in the particular case where~$u(t)\equiv \text{TRUE}$
the solution to the minimal observability problem is to make~$8$ SV directly measurable, namely,~$X_1 ,\dots,X_8 $. 
Note that using the STP is applicable here, as there are only nine~SVs
yet the transition matrix is already~$2^9\times 2^9$. 

Applying Algorithm~2 to the BCN~\eqref{eq:boolean_model}  (without assuming a necessarily 
constant control) 
 generates the lists:
$L_1=\{X_1,\dots,X_8\}$,~$L_2=\{X_9\}$, and~$L_C=\emptyset$, and then returns~$L_1$. 
Thus, in this particular
 case Algorithm~2 provides an optimal solution. We emphasize again that other algorithms, 
that require exponential run-time, are not applicable for~BCNs with, say, $n\geq 30$ SVs.

 \end{document}